\documentclass[aap,MSNbibl,dvips]{arximspdf}
\usepackage{mathrsfs}

\doi{10.1214/13-AAP928} 
\volume{24}
\issue{1}
\pubyear{2014}
\firstpage{447}
\lastpage{448}

\makeatletter

\renewcommand{\P}{\mathbb{P}}

\makeatother

\begin{document}
\begin{frontmatter}

\title{A note on ``Bayesian nonparametric estimators derived from
conditional Gibbs structures''}
\runtitle{Note}

\begin{aug}
\author[A]{\fnms{Antonio} \snm{Lijoi}\corref{}\ead[label=e1]{lijoi@unipv.it}},
\author[B]{\fnms{Igor} \snm{Pr\"unster}\ead[label=e2]{igor.pruenster@unito.it}}
\and
\author[C]{\fnms{Stephen G.} \snm{Walker}\ead[label=e3]{s.g.walker@math.utexas.edu}}
\runauthor{A. Lijoi, I. Pr\"unster and S. G. Walker}
\affiliation{University of Pavia and Collegio Carlo Alberto,
Moncalieri,\break
University of Torino and Collegio Carlo Alberto, Moncalieri,
and~University of Texas at Austin}
\address[A]{A. Lijoi\\
Department of Economics \\
\quad and Management\\
University of Pavia\\
Via San Felice 5\\
27100 Pavia\\
Italy \\
\printead{e1}}
\address[B]{I. Pr\"unster\\
Department of Economics\\
\quad and Statistics\\
University of Torino\\
Corso Unione Sovietica 218/bis\\
I-10134 Torino\\
Italy \\
\printead{e2}}
\address[C]{S. G. Walker\\
Department of Mathematics\\
University of Texas at Austin\\
1 University Station C1200\\
Austin, Texas 78712\\
USA\\
\printead{e3}}
\end{aug}

\received{\smonth{2} \syear{2013}}



\end{frontmatter}

The present note aims at clarifying some possibly confusing notation
used in~\cite{Lij08}, even if its correct interpretation should be
clear from context and from the proofs contained therein. In
particular, the expression 
%
\renewcommand{\theequation}{$\ast$}
%
\begin{equation}\label{eq2}
\P\bigl[K_n=k,\mathbf{N}_{n}=(n_1,
\ldots,n_{K_n})\bigr]
\end{equation}
displayed in (3) of \cite{Lij08} has been used to indicate the
probability of observing a specific realization of a random partition
of the integers $[n]=\{1,\ldots,n\}$ into $k$ blocks of sizes
$(n_1,\ldots,n_k)\in\Delta_{n,k}$. However, the notation in
(\ref{eq2}) may be actually misleading since its natural
interpretation is as the probability of all partitions of
$[n]=\{1,\ldots,n\}$ into $k$ blocks of sizes
$(n_1,\ldots,n_k)\in\Delta_{n,k}$. For this reason, (3) in
\cite{Lij08} should be rewritten as
\renewcommand{\theequation}{$3$}
%
\begin{equation}\label{eq3}
\Pi_k^{(n)}(n_1, \ldots,n_k)=\frac{\theta^k}{(\theta)_n} \prod_{j=1}^k
(n_j-1)!,
\end{equation}
where $\Pi_k^{(n)}$ is the \textit{exchangeable partition probability
function} notation introduced in Section~2. Hence, if $\Pi_n$ is a
random element taking values in the set of all partitions of $[n]$, one
has $\Pi_k^{(n)}(n_1,\ldots,n_k)=\P[\Pi_n=\pi]$ for any partition $\pi$
of $[n]$ into $k$ blocks $\{A_1,\ldots,A_k\}$ with $|A_i|=n_i$, for
$i=1,\ldots,k$. See also the first displayed equation on page 1522
after (4) in \cite{Lij08}. The probability of all partitions of $[n]$
into $k$ blocks with respective sizes $n_1,\ldots,n_k$ (in exchangeable
order) is then equal to $n! \Pi_k^{(n)}(n_1,\ldots,n_k)/(k!n_1! \cdots
n_k!)$. The same caveat also applies to (28) and (44)--(46) in
\cite{Lij08}.

An analogous clarification concerns
%
\renewcommand{\theequation}{$\ast\ast$}
%
\begin{equation}\label{eq1}
\P\bigl[K_m^{(n)}=k, L_m^{(n)}=s,
\mathbf{S}_{L_m^{(n)}}=(s_1,\ldots,s_{K_m^{(n)}}) | K_n=j\bigr]
\end{equation}
displayed in (9) of \cite{Lij08}. Indeed, (\ref{eq1}) has been used
to denote the probability of a specific partition of $[s]$ into $k$
blocks with sizes $(s_1,\ldots,s_k)$ in $\Delta_{s,k}$, given any
partition of $[n]$ into $j$ parts, thus differing from its natural
meaning as the probability of all such partitions. Our intepretation of
(\ref{eq1}) in \cite{Lij08} is clearly consistent with the
right-hand side of (9), that is,
\[
\frac{V_{n+m,j+k}}{V_{n,j}} \pmatrix{m
\cr
s} (n-j\sigma)_{m-s} \prod
_{i=1}^k(1-\sigma)_{s_i-1}
\]
and with the result displayed in Corollary 1, where equation (10)
arises after multiplying the right-hand side of (9) by $s!/(k!s_1!
\cdots s_k!)$ and, then, marginalizing with respect to all
$(s_1,\ldots,s_k)$ in $\Delta_{s,k}$. Similar remarks apply to the notation
appearing in the following displayed formulas: (17), (19), equation at
the end of page 1528, (21), equation right after Corollary 2 on page
1529, (22), equation at the end of page 1531, (34).




\printaddresses

\end{document}